\documentclass{amsart}
\usepackage{amsthm,amsmath,amsfonts}
\usepackage{mathrsfs}
\usepackage[colorlinks=false,linktocpage=true]{hyperref}
\usepackage{tocvsec2}
\usepackage{latexsym}
\numberwithin{equation}{section}
\font\tenscrpt=eusm10
\font\sevenscrpt=eusm10 scaled 700
\font\fivescrpt=eusm10 scaled 500
\newfam\eusmfam
\textfont\eusmfam=\tenscrpt
\scriptfont\eusmfam=\sevenscrpt
\scriptscriptfont\eusmfam=\fivescrpt
\def\scr#1{{\fam\eusmfam\relax#1}}
\def\qed{\quad\vcenter{\hrule\hbox{\vrule height.6em\kern.6em\vrule}\hrule}}

\newenvironment{pf*}[1]{{\textsc #1.\quad}}{$\qed$\bigskip\newline}
\theoremstyle{plain}
\newtheorem{thm}{Theorem}[section]

\newtheorem{defn}{Definition}[section]
\newtheorem{rem}{Remark}[section]
\newcommand{\thmref}[1]{Theorem~\ref{#1}}

\newcommand{\remref}[1]{Remark~\ref{#1}}
\newcommand{\eqnref}[1]{{\rm (\ref{#1})}}

\def\square{\quad\vcenter{\hrule\hbox{\vrule height.6em\kern.6em\vrule}\hrule}}
\numberwithin{equation}{section}
\font\tenscrpt=eusm10
\font\sevenscrpt=eusm10 scaled 700
\font\fivescrpt=eusm10 scaled 500
\newfam\eusmfam
\textfont\eusmfam=\tenscrpt
\scriptfont\eusmfam=\sevenscrpt
\scriptscriptfont\eusmfam=\fivescrpt
\def\scr#1{{\fam\eusmfam\relax#1}}
\def\eqdef{\overset{\triangle}{=}}
\def\ehNabh{e_{heat}^{Neu}(a,b,h)}
\def\ehNabdh{e_{heat}^{Neu}(a,b+d,h)}
\def\ehNazh{e_{heat}^{Neu}(a,0,h)}
\def\RNbTtauniLYiWit{\Upsilon _{T\wedge\tauni}^{\RVi,\tilde {\Wm}^{(i)}}([0,L])}
\def\RNbttauniBYiWit{\Upsilon _{t\wedge\tauni}^{\RVi,\tilde {\Wm}^{(i)}}(B)}

\def\RNfTtauniLYiWi{\Xi_{T\wedge\tauni}^{\RVi,{\Wm}^{(i)}}([0,L])}

\def\RNfTtaunoLYoWo{\Xi _{T\wedge\tauno}^{\RVo,{\Wm}^{(1)}}([0,L])}
\def\RNfTtauntLYtWt{\Xi _{T\wedge\taunt}^{\RVt,{\Wm}^{(2)}}([0,L])}

\def\N{{\Bbb N}}
\def\R{{\Bbb R}}
\def\Rpt{{\Bbb R}_+^2}

\def\L{{\Bbb L}}

\def\P{{\Bbb P}}
\def\Pti{\tilde{\P}}
\def\Pii{\P^{(i)}}
\def\Pni{\P_n^{(i)}}
\def\Piti{\tilde{\P}^{(i)}}
\def\Poti{\tilde{\P}^{(1)}}
\def\Ptti{\tilde{\P}^{(2)}}
\def\Pno{\P_n^{(1)}}
\def\Pnt{\P_n^{(2)}}
\def\EP{{\Bbb E}_{\P}}
\def\EPno{{\Bbb E}_{\Pno}}
\def\EPnt{{\Bbb E}_{\Pnt}}
\def\Ft{{\scr{F}}_t}

\def\FTi{{\scr{F}}_T^{(i)}}
\def\filpspace{(\Omega, \scr{F}, \{{\scr{F}}_t\},\P)}
\def\filpspacei{(\Omega^{(i)}, {\scr{F}}^{(i)}, \{{\scr{F}}_t^{(i)}\},\Pii)}

\def\filpspaceTni{(\Omega^{(i)}, \scr{F}_T^{(i)}, \{{\scr{F}}_t^{(i)}\},{\P}_n^{(i)})}

\def\filpspacetiti{(\tilde{\Omega}, \tilde{\scr{F}}, \{\tilde{{\scr{F}}}_t\},\Pti)}
\def\filpspaceiti{(\Omega^{(i)}, {\scr{F}}^{(i)}, \{{\scr{F}}_t^{(i)}\},\Piti)}
\def\Ru{R_u}
\def\RU{R_U}
\def\RV{R_V}
\def\RVi{R_{\Vi}}
\def\RVo{R_{\Vo}}
\def\RVt{R_{\Vt}}
\def\RZ{R_Z}

\def\Vo{V^{(1)}}

\def\Vt{V^{(2)}}

\def\Ui{U^{(i)}}
\def\Vi{V^{(i)}}

\def\Wm{\scr W}
\def\Wim{\Wm^{(i)}}

\def\Wtim{\tilde{\Wm}}
\def\Witim{\Wtim^{(i)}}
\def\solViWin{(\Vi,\Wmin)}

\def\solViWiti{(\Vi,{\Witim})}
\def\LawUiPi{\L_{\Pii}^{\Ui}}
\def\tauni{\tau_n^{(i)}}
\def\tauno{\tau_n^{(1)}}
\def\taunt{\tau_n^{(2)}}

\def\Wmi{\Wm^{(i)}}
\def\Wmit{\tilde{\Wm}^{(i)}}
\def\Wmin{\Wm_n^{(i)}}
\def\WitB{\Wm_t^{(i)}(B)}
\def\WitBn{\Wm_{t\wedge\tauni}^{(i)}(B)}

\def\WittB{\tilde{\Wm}_t^{(i)}(B)}
\def\BL{{\scr B}([0,L])}
\def\vfi{\varphi}
\def\vfidp{\vfi''}
\def\vfix{\vfi(x)}
\def\RtauniTL{{\scr R}_{T\wedge\tauni,L}}
\def\RtL{{\scr R}_{t,L}}
\def\RTL{{\scr R}_{T,L}}

\def\RTLTcl{{\overset{\smile}{\scr R}}_{T,L}}
\def\CRTLR{C(\RTL,\R)}

\begin{document}
\title[Uniqueness in law for the Allen-Cahn SPDE]{Uniqueness in Law for the Allen-Cahn SPDE via Change of Measure}
\author {Hassan Allouba}
\address{Department of Mathematics and Statistics, University of Massachusetts,
Amherst, MA 01003-4515 \\
E-mail: allouba@math.umass.edu, Phone: (413) 545-6010, Fax: (413) 545-1801}
\date{January 2000}
\begin{abstract}
We start by first using change of measure to prove the transfer of
uniqueness in law among pairs of parabolic SPDEs differing only by a drift function,
under an almost sure $L^2$ condition
on the drift/diffusion ratio.  This is a considerably weaker condition than the
usual Novikov one, and it allows us to prove
uniqueness in law for the Allen-Cahn SPDE driven by space-time
white noise with diffusion function $a(t,x,u)=Cu^\gamma$, $1/2\le\gamma\le1$ and $C\ne0$.
The same transfer result is also valid for ordinary SDEs
and hyperbolic SPDEs.
\end{abstract}
\maketitle
\section{Introduction.}
We start by considering the pair of parabolic SPDEs 
\begin{equation}\label{Pb}
\begin{cases}
\displaystyle\frac{\partial U}{\partial t}=\Delta _{x}U+b(t,x,U)+a(t,x,U)
\displaystyle\frac{\partial^2 W}{\partial t\partial x};
& (t,x)\in\RTLTcl,
\cr U_x(t,0)=U_x(t,L)=0; & 0<t\le T,
\cr U(0,x)=h(x); & 0<x<L,
\end{cases}
\end{equation}
and
\begin{equation}
\begin{cases}
\displaystyle\frac{\partial V}{\partial t}=\Delta _{x}V+(b+d)(t,x,V)+a(t,x,V)
\displaystyle\frac{\partial^2 W}{\partial t\partial x};
& (t,x)\in\RTLTcl,
\cr V_x(t,0)=V_x(t,L)=0; & 0<t\le T,
\cr V(0,x)=h(x); & 0<x<L,
\end{cases}
\label{Pbd}
\end{equation}
on the space-time rectangle $\RTL\eqdef[0,T]\times[0,L]$, where
$\RTLTcl\eqdef(0,T]\times(0,L)$.  $W(t,x)$ is
the Brownian sheet corresponding to the driving space-time white noise,
written formally as $\partial^2W/\partial t\partial x$.  As in Walsh \cite{WA},
white noise is regarded as a continuous orthogonal martingale
measure, which we denote by $\Wm$, with the corresponding Brownian sheet
as the random field induced by $\Wm$ in the usual way:
$W(t,x)=\Wm([0,t]\times[0,x])\eqdef \Wm_t([0,x])$.  The diffusion
$a(t,x,u)$ and the drifts $b(t,x,u)$ and $d(t,x,u)$ are
Borel-measurable $\R$-valued functions on $\RTL\times\R$; and
$h:\RTL\to\R$ is a Borel-measurable function.
Henceforth, we will denote \eqnref{Pb} and
\eqnref{Pbd} by $\ehNabh$ and $\ehNabdh$, respectively.  When $b\equiv0$, we denote
\eqnref{Pb} by $\ehNazh$.  In the interest of getting quickly to our main results,
we relegate to the Appendix the rigorous interpretation of all SPDEs
considered in this paper.

Proceeding toward a precise statement of the main result,
we adopt some convenient notation.  Let $\Ru(t,x)\eqdef d(t,x,u)/a(t,x,u)$,
for any $(t,x,u)\in[0,T]\times[0,L]\times\R$, whenever the ratio is well defined.
Let $\lambda$ denote Lebesgue measure on $\Rpt$.  Our main result for
the pair $\ehNabh$ and $\ehNabdh$ can now be stated as
\begin{thm}
Assume that $\RU$ and $\RV$ are in $L^2(\RTL,\lambda)$, almost surely,
whenever the random fields $U$ and $V$ solve $($weakly or strongly\/$)$
$\ehNabh$ and $\ehNabdh$, respectively $($see \remref{solab}\/$)$.  Then,
uniqueness in law holds for $\ehNabh$ iff
uniqueness in law holds for $\ehNabdh$.
\label{Partrnsfr}
\end{thm}
By transferring uniqueness in law under weaker conditions
(almost sure $L^2$ vs. Novikov's),
\thmref{Partrnsfr} makes more applicable the
notion of Girsanov equivalence in our earlier work
(Theorem 3.3.2 in \cite{HAd} or Theorem 4.2 in \cite{HA1}).  The Neumann conditions in $\ehNabh$ and
$\ehNabdh$ may be changed to Dirichlet
conditions without affecting the conclusions of \thmref{Partrnsfr}.

An interesting application of \thmref{Partrnsfr} is provided in
\thmref{ACstab} below
for the stochastic Allen-Cahn equation driven by space-time white noise
\begin{equation}
\begin{cases}
\displaystyle\frac{\partial V}{\partial t}=\Delta _{x}V+2V(1-{V}^{2})+ C V^\gamma
\displaystyle\frac{\partial^2 W}{\partial t\partial x};
& (t,x)\in\RTLTcl,
\cr V_x(t,0)=V_x(t,L)=0; & 0<t\le T,
\cr V(0,x)=h(x); & 0<x<L,
\end{cases}
\label{AC}
\end{equation}
 in the case $C\ne0$ and $\frac12\le\gamma\le1$.  The deterministic
Allen-Cahn PDE was introduced by Allen and Cahn \cite{AllenC} as a model for grain
boundary motion.  It has since become an important PDE for many
mathematicians (see e.g. Katsoulakis et al.~\cite{Markos}, Sowers et al.~\cite{Sowers}, and
the references therein); and we intend to investigate, in a future paper, further properties of its
solutions in the presence of a driving space-time white noise.
We are thankful to Markos Katsoulakis for interesting Allen-Cahn conversations.
\begin{thm}
Consider the stochastic Allen-Cahn equation $\eqnref{AC}$ $(C\neq0)$.
If $\gamma\in[\frac12,1]$ then uniqueness in law holds for $\eqnref{AC}$.
\label{ACstab}
\end{thm}
\thmref{ACstab}  follows since
\begin{itemize}
\item The Allen-Cahn SPDE \eqnref{AC} satisfies our transfer condition, and
\item A result of Mytnik \cite{L} gives us uniqueness in law for the SPDE \eqnref{Pb} with $b\equiv0$ and $a(t,x,u)=Cu^\gamma$, $\gamma\in(1/2,1)$, the case $a(t,x,u)=Cu^{1/2}$ admits uniqueness in law as discussed in \cite{Mu} p.~326 and in
\cite{RC}, and a classic result of Walsh gives us uniqueness for \eqnref{Pb} when $b\equiv0$ and $a(t,x,u)=Cu$.
\end{itemize}
With minor adaptations, the same uniqueness transfer result
in \thmref{Partrnsfr} holds for ordinary SDEs and
hyperbolic (wave) SPDEs (see \cite{HA1} Theorem 3.6, Theorem 5.2, and their proofs
for our uniqueness and existence transfer result for space-time SDEs and
their rotationally-equivalent wave SPDEs, using change of measure under
Novikov's condition).
We note that the existence of solutions to heat SPDEs with continuous
diffusion coefficients $a$ satisfying a linear growth condition was established in
\cite{HA2,Mu,Re}.  In \cite{HA2}, we used an approximating system of stochastic
differential-difference equations (SDDEs) to give a new proof of Reimers'
existence result, then we used our Girsanov theorem from \cite{HA1} to
extend the result to measurable drifts, under Novikov's condition.
\section{Proof of the Main Result.}  We begin by adapting the well known
Novikov condition to our setting:
we say that a predictable random field $X$ on the probability space
$\filpspace$ (see Walsh \cite{WA}) satisfies Novikov's condition on $\RTL$ if
\begin{equation}
\EP\left[\exp\left(\frac12\intop _{\RTL}X^2(t,x)dtdx\right)\right]<\infty.
\label{N}
\end{equation}
\begin{rem}
It is clear that if
$\Ru(t,x)$ is uniformly bounded for $(t,x,u)\in\RTL\times\R$, then
$\RZ$ satisfies Novikov's condition on $\RTL$ for every predictable random field $Z$.
\label{bdd}
\end{rem}

{\sc Proof of \thmref{Partrnsfr}}
Assume that uniqueness in law holds for $\ehNabh$, and suppose that
$$\solViWiti,\ \filpspaceiti;\ i=1,2,$$
are solutions to $\ehNabdh$.  By assumption
\begin{equation}
\Piti\Bigg[\int _{\RTL}{R_{V^{(i)}}^{2}}(t,x)dtdx<\infty \Bigg]=1,\ i=1,2.
\label{L2}
\end{equation}
Now take $\{\tauni\}$ to be the sequence of stopping times
\begin{equation}
\tauni\eqdef T\wedge\inf\left\{0\le t\le T;\intop_{\RtL}\RVi^2(s,x)ds dx=n\right\};
\ n\in\N,\ i=1,2.
\label{STimes}
\end{equation}
Let $\Wmi=\{\WitB,\Ft;0\le t\le T,B\in\BL\}$
be given by
$$\WitB\eqdef\WittB+\int _{[0,t]\times B}\RVi(s,x)dsdx;\ i=1,2.$$
Novikov's condition \eqnref{N} and Girsanov's theorem for white noise (see
Corollary 3.1.3 in \cite{HAd}) imply that $\Wmin=\{\WitBn,\Ft;0\le t\le T,B\in\BL\}$
is a white noise stopped at time $\tauni$, under the probability measure
$\Pni$ defined on $\FTi$ by the recipe
\begin{equation*}
\frac{d\Pni}{d\Piti}=\RNbTtauniLYiWit;\ n\in\N,\ i=1,2,
\end{equation*}
where
\begin{equation*}
\begin{split}
&\RNbttauniBYiWit
\\&\eqdef\ \exp\left[-\intop_{[0,t\wedge\tauni]\times B}\RVi(s,x)\right.\Wmit(ds,dx)
\left.- \frac12\intop_{[0,t\wedge\tauni]\times B} \RVi^2(s,x) ds dx\right];
\end{split}
\end{equation*}
$0\le t\le T$, $B\in\BL$.  It follows that $\solViWin$, $\filpspaceTni$ is a solution
to $\ehNabh$ on $\RtauniTL\eqdef[0,T\wedge\tauni ]\times[0,L]$ for each $i=1,2$ and $n\in\N$.
Of course, for $i=1,2$,
\begin{equation}
\begin{split}
&\frac{d\Piti}{d\Pni}=\RNfTtauniLYiWi\\
\eqdef&\exp\left[\intop_{[0,T\wedge\tauni]\times[0,L]}\RVi(s,x)\right.\Wmi(ds,dx)
\left.- \frac12\intop_{[0,T\wedge\tauni]\times[0,L]} \RVi^2(s,x) ds dx\right];
\end{split}
\label{RDF}
\end{equation}
$n\in\N$.  Consequently, for any set $\Lambda\in{\scr B}(\CRTLR)$
\begin{equation}
\begin{split}
\Poti\left[\Vo\in\Lambda,\tauno=T\right]
&=\EPno\left[1_{\{\Vo\in\Lambda,\tauno=T\}}{\RNfTtaunoLYoWo}\right]\\
&=\EPnt\left[1_{\{\Vt\in\Lambda,\taunt=T\}}{\RNfTtauntLYtWt}\right]\\
&=\Ptti\left[\Vt\in\Lambda,\taunt=T\right];\ \forall n\in\N,
\end{split}
\label{uil}
\end{equation}
where we have used the uniqueness in law assumption on $\ehNabh$ (comparing
the $\Vi$'s only on $\Omega_n^{(i)}\eqdef\{\tauni=T\}$ for each $n$), \eqnref{STimes},
and \eqnref{RDF} to get the second equality in \eqnref{uil}.
By \eqnref{L2} and \eqnref{STimes}, we get that
$\lim \nolimits _{n\rightarrow \infty }\Piti[\tauni=T]=1\mbox{ for } i=1,2$.
We then see that passing to the limit as $n\rightarrow \infty $ in \eqnref{uil} gives us that the
law of $\Vo$ under $\Poti$ is the same as that of $\Vt$ under $\Ptti$.  I.e.,
we have uniqueness in law for $\ehNabdh$.  The proof of the other direction
is similar and is omitted.
$\square$

Our Uniqueness result for the Allen-Cahn SPDE \eqnref{AC}
can now be proved.
\vspace{2mm} \\
{\sc Proof of \thmref{ACstab}}
By \thmref{Partrnsfr}, the proof essentially reduces to checking
whether the random fields $\RU$ and $\RV$
are in $L^2(\RTL,\lambda)$, almost surely,
whenever $U$ solves (weakly or strongly)
$\ehNazh$ (with $a(t,x,u)\equiv C u^\gamma$ and $\frac12\le\gamma\le1$)
and $V$ solves (weakly or strongly) the Allen-Cahn SPDE \eqnref{AC}.
That this is true can easily be seen since, in this case,
\begin{equation}
\RU^2(t,x)=\frac{4}{C^2}U^{2(1-\gamma)}(U^4-2U^2+1),\mbox{ and }
\RV^2(t,x)=\frac{4}{C^2}V^{2(1-\gamma)}(V^4-2V^2+1).
\label{XY}
\end{equation}
The continuity of $U$ and $V$ implies that $\RU^2$ and $\RV^2$ are continuous, for any $0\le\gamma\le1$.
Therefore,
if $U$ and $V$ are defined on the usual probability spaces $\filpspace$ and $\filpspacetiti$,
respectively, then
\begin{equation*}
\begin{split}
|\RU^2(t,x,\omega)|\le K(\omega,\gamma)<\infty;\mbox{ for all } (t,x)\in\RTL,\ \gamma\in[0,1]\mbox{ a.s. }\P,\\
|\RV^2(t,x,\tilde\omega)|\le\tilde {K}(\tilde\omega,\gamma)<\infty;\mbox{ for all } (t,x)\in\RTL,\
\gamma\in[0,1] \mbox{ a.s. }\Pti,
\end{split}
\end{equation*}
where $K$ and $\tilde K$ depend only on $(\omega,\gamma)\in\Omega\times[0,1]$ and
$(\tilde\omega,\gamma)\in\tilde\Omega\times[0,1]$, respectively.  It follows that, for any fixed but arbitrary
$\gamma\in[0,1]$, $\RU$ and $\RV$ are in
$L^2(\RTL,\lambda)$, almost surely.  The assertion of \thmref{ACstab} then follows from
\thmref{Partrnsfr} and the fact that
uniqueness in law holds for $\ehNazh$ when $a(t,x,u)\equiv C u^\gamma$ and
$\gamma\in[\frac12,1]$ (see \cite{Mu,RC,L} and \cite{WA}).
$\square$
\appendix
\message{Appendix A}
\setcounter{section}{1}
\section*{Appendix}
We collect here definitions and conventions that are used throughout this
article.
Filtrations are assumed to satisfy the usual conditions
(completeness and right continuity), and any
probability space $\filpspace$ with such a filtration is termed a usual
probability space.  The space of continuous functions on $\RTL$ is denoted by $C(\RTL)$.
\begin{defn}[Strong and Weak Solutions to $\ehNabh$]
We say that the pair $(U,\Wm)$ defined on the usual probability space
$\filpspace$ is a solution to the stochastic heat equation $\ehNabh$ if
$\Wm$ is a space-time white noise on $\RTL$;
the random field $U(t,x)$ is predictable $($as in \cite{WA}\/$)$, with
continuous paths on $\RTL$; and
the pair $(U,\Wm)$ satisfies the test function formulation:
\begin{equation*}
\begin{split}
&\int_0^L (U(t,x)-h(x))\vfix dx-\int_0^L\int_0^t U(s,x)\vfidp(x)ds dx
\\=&\int_0^L\int_0^t \left[a(s,x,U(s,x))\vfix \Wm(ds,dx)+
 b(s,x,U(s,x))\vfix ds dx\right];\ 0\leq t\leq T,
\end{split}
\label{TFF}
\end{equation*}
$\mbox{ a.s. } \P$,  for every
$\vfi\in\Theta _{0}^{L}\eqdef\left\{ \varphi \in C_c^\infty (\R;\R):
\varphi \prime (0)=\varphi \prime (L)=0\right\}$ $($$C_c^\infty (\R;\R)$ being
the collection of smooth $\R$-valued function on $\R$ with compact support$)$.
A solution
is said to be strong if the white noise $\Wm$ and  the usual probability
space
$\filpspace$ are fixed a priori and $\Ft$ is the augmentation of the natural
filtration for $\Wm$ under $\P$.  It
is termed a weak solution if we are allowed to choose the usual
probability space and the white noise $\Wm$ on it, without requiring that
the filtration be the augmented natural filtration of $\Wm$.
\label{heatsol}
\end{defn}
\begin{rem}
We often simply say that $U$ solves $\ehNabh$ (weakly or strongly) to mean the
same thing as above.
\label{solab}
\end{rem}
\begin{defn}[Uniqueness for SPDEs] We say
that uniqueness in law holds for $\ehNabh$ if the laws $\LawUiPi$ of $\Ui$
under $\Pii$; $i=1,2$, are the same on $(C(\RTL),{\scr B}(C(\RTL))$ whenever
$(\Ui,\Wim)$, $\filpspacei$; $i=1,2$, are solutions to $\ehNabh$.
\label{pathlawuniqSPDES}
\end{defn}
 
{\bf Acknowledgements.}  The author would like to thank the referee for her comments which improved the presentation
of this paper.

\end{document}